\documentclass[aap,MSNbibl,citesort,dvips]{arximspdf}

\doi{10.1214/10-AAP715}
\volume{21}
\issue{3}
\pubyear{2011}
\firstpage{987}
\lastpage{1016}

\makeatletter

\newtheorem{theorem}{Theorem}[section]
\newtheorem{lemma}[theorem]{Lemma}

\newproclaim{definition}[theorem]{Definition}
\newproclaim{remark}[theorem]{Remark}

\newcommand{\bigtimes}{\mathop{\mbox{\fontsize{17}{17}\selectfont{$\!\times$}}}}

\newcommand{\E}{\mathcal{E}}
\newcommand{\N}{\mathcal{N}}
\newcommand{\G}{\mathcal{G}}
\newcommand{\I}{\mathcal{I}}
\newcommand{\eps}{\varepsilon}
\newcommand{\BR}{\operatorname{BR}}
\newcommand{\F}{\mathcal{F}}
\newcommand{\D}{\mathcal{D}}

\newcommand{\Prob}{\mathbb{P}}

\newcommand{\PrG}{\mathbb{P}_\mathcal{G}}
\newcommand{\PrGT}{\mathbb{P}_{(n,p)}}
\newcommand{\PrGG}{\mathbb{P}_{G}}
\newcommand{\EGG}{\mathbb{E}_{G}}
\newcommand{\ones}{\mathbb{I}}

\newcommand{\stexp}{\mathbb{E}}
\newcommand{\M}{\mathcal{M}}

\makeatother

\begin{document}
\begin{frontmatter}

\title{Connectivity and equilibrium in random games}
\runtitle{Connectivity and equilibrium in random games}

\begin{aug}
\author[A]{\fnms{Constantinos} \snm{Daskalakis}\thanksref{t1}\ead[label=e1]{costis@csail.mit.edu}},
\author[B]{\fnms{Alexandros G.} \snm{Dimakis}\corref{}\thanksref{t2}\ead[label=e2]{dimakis@usc.edu}}\\ and
\author[C]{\fnms{Elchanan} \snm{Mossel}\thanksref{t3}\ead[label=e3]{mossel@stat.berkeley.edu}}
\runauthor{C. Daskalakis, A. G. Dimakis and E. Mossel}
\affiliation{Massachusetts Institute of Technology, University of
Southern California, Los~Angeles~and University of California, Berkeley}
\address[A]{C. Daskalakis\\
Computer Science and Artificial\\
\quad Intelligence Laboratory\\
Department of Electrical Engineering \\
\quad and Computer Science\\
Massachusetts Institute of Technology\\
Cambridge, Massachusetts 02139\\
USA\\
\printead{e1}}
\address[B]{A. G. Dimakis\\
Department of Electrical Engineering \\
University of Southern California\\
Los Angeles, California 90089\\
USA\\
\printead{e2}}
\address[C]{E. Mossel\\
Department of Statistics \\
University of California, Berkeley\\
Berkeley, California 94720\\
USA\\
\printead{e3}}
\end{aug}

\thankstext{t1}{Supported by a Sloan Fellowship in Computer Science and
NSF Grant CCF-0953960 (CAREER).
Part of this work was done while the author was at EECS, UC Berkeley,
with the support of a Microsoft Research Fellowship and NSF Grants
CCF-0635319, DMS-05-28488 and DMS-05-48249 (CAREER).}
\thankstext{t2}{Part of this work was done while the author was at
EECS, UC Berkeley. Supported by NSF Grant DMS-05-28488 and a Microsoft
Research Fellowship.}
\thankstext{t3}{Supported by a Sloan fellowship in Mathematics, NSF
Grants DMS-05-28488 and DMS-05-48249 (CAREER) and ONR Grant N0014-07-1-05-06.}

\received{\smonth{12} \syear{2007}}
\revised{\smonth{5} \syear{2010}}

%
\begin{abstract}
We study \textit{how the structure of the interaction graph} of a game
affects the existence of pure Nash equilibria. In particular, for a
fixed interaction graph, we are interested in whether there are pure
Nash equilibria arising when random utility tables are assigned to the
players. We provide conditions for the structure of the graph under
which equilibria are likely to exist and complementary conditions which
make the existence of equilibria highly unlikely. Our results have
immediate implications for many deterministic graphs and generalize
known results for random games on the complete graph. In particular,
our results imply that the probability that bounded degree graphs have
pure Nash equilibria is exponentially small in the size of the graph
and yield a simple algorithm that finds small nonexistence
certificates for a large family of graphs. Then we show that in any
strongly connected graph of $n$ vertices with expansion $(1+\Omega(1))
\log_2 (n)$ the distribution of the number of equilibria approaches the
Poisson distribution with parameter $1$, asymptotically as $n \to+
\infty$.

In order to obtain a refined characterization of the degree of
connectivity associated with the existence of equilibria, we also study
the model in the random graph setting. In particular, we look at the
case where the interaction graph is drawn from the Erd\H{o}s--R\'enyi,
$G(n,p)$, model where each edge is present independently with
probability $p$. For this model we establish a \textit{double phase
transition} for the existence of pure Nash equilibria as a function of
the average degree $p n$, consistent with the nonmonotone behavior of
the model. We show that when the average degree satisfies $n p > (2 +
\Omega(1)) \log_e (n)$, the number of pure Nash equilibria follows a
Poisson distribution with parameter~$1$, asymptotically as $n
\rightarrow\infty$. When $1/n \ll n p < (0.5 -\Omega(1)) \log_e (n)$,
pure Nash equilibria fail to exist with high probability. Finally, when
$n p = O(1/n)$ a pure Nash equilibrium exists with constant
probability.
%
%
\end{abstract}

%

\end{frontmatter}

\section{Introduction}


In recent years there has been a convergence of ideas from computer
science and the social sciences aiming to model and analyze large
complex networks such as the web graph, social networks and
recommendation systems. From the computational perspective, it has been
recognized that the successful design of algorithms performed on such
networks, including routing, ranking and recommendation algorithms,
must take into account the social dynamics and economic incentives as
well as the technical properties that govern these networks
\cite
{DBLPconfstocPapadimitriou01,DBLPconfstocRaghavan06,DBLPconffocsKleinberg06}.

Game theory has been very successful in modeling strategic behavior in
large systems of economically incentivized entities. In the context of
routing, for instance, it has been employed to study the effect of
selfishness on the efficiency of a network, whereby the performance of
the network at equilibrium is compared to the performance when a
central authority can simply dictate a solution
\cite
{roughgarden02price,roughgardenTardos2002,roughgardenTardos2003,correa04selfish}.
The effect of selfishness has
been studied in several other settings, for example, load balancing
\cite{czumaj02selfish,545436,koutsoupias99worstcase,roughgarden01stackelberg},
facility location~\cite{vetta02nash} and
network design~\cite{anshelevich04price}.

A simple way to model interactions between agents in a large network is
with a \textit{graphical game}~\cite{kearns}: a graph $G=(V,E)$ is
defined whose vertices represent the players of the game, and an edge
$(v,w)\in E$ corresponds to the strategic interaction between players
$v$ and $w$; each player $v \in V$ has a finite set of strategies
$S_v$, which throughout this paper will be assumed to be binary so that
there are two possible strategies for each player.
A \textit{utility}, \textit{or payoff}, \textit{table} $u_v$ for player
$v$ assigns a real
number $u_v(\sigma_v, \sigma_{\N(v)})$ to every selection of
strategies by player $v$ and the players in $v$'s neighborhood, that
is, the set of nodes $v'$ such that $(v,v') \in E$, denoted by $\N
(v)$. A~\textit{pure Nash equilibrium} (PNE) of the game is some state,
or \textit{strategy profile}, $\sigma$ of the game, assigning to every
player $v$ a single strategy $\sigma_v \in S_v$, such that no player
has a unilateral incentive to deviate. Equivalently, for every player
$v \in V$,
%
%
\begin{equation}\label{eq: equilibrium constraint}
u_v\bigl(\sigma_v,\sigma_{\N(v)}\bigr) \ge u_v\bigl(\sigma
'_v,\sigma
_{\N(v)}\bigr)\qquad\mbox{for every strategy $\sigma_v' \in S_v$}.
\end{equation}
When condition (\ref{eq: equilibrium constraint}) is satisfied, we say
that the strategy $\sigma_v$ is a \textit{best response to the strategies
$\sigma_{\N(v)}$}.

The concept of the pure strategy Nash equilibrium is more compelling,
decision theoretically,
than the concept of the mixed strategy Nash equilibrium, its
counterpart that allows players to choose distributions over their
strategy sets. This is because it is not always meaningful in
applications to assume that the players of a game may adopt randomized
strategies.
Unfortunately, unlike mixed Nash equilibria, PNE do not always exist.
It is then an important problem to study how the existence of PNE
depends on the properties of the game.

The focus of this paper is to understand how the connectivity of the
underlying graph affects the existence of a PNE. We obtain two kinds of
results. The first concerns the existence of a PNE in an ensemble of
random graphical games defined on a random, $G(n,p)$, graph. We obtain
a characterization of the probability that a PNE exists as a function
of the density of the graph. The second set of results concerns random
graphical games on deterministic graphs. Here, we obtain conditions on
the structure of the graph under which a PNE does not exist with high
probability, suggesting also an efficient algorithm for finding
witnesses of the nonexistence of a PNE. We also give complementary
conditions on the structure of the graph under which a PNE exists with
constant probability. Our results are described in detail in
Section~\ref{sec: main results}.

\subsection*{Comparison to typical constraint satisfaction problems}
Graphical games provide a more compact way of representing large
networks of interacting agents than normal form games, in which the
game is described as if it were played on the complete graph. Besides
the compact description, one of the motivations
for the Introduction of graphical games is their intuitive affinity to
graphical statistical models; indeed, several algorithms for graphical
games do have the flavor of algorithms for solving Bayes nets or
constraint satisfaction problems \mbox{\cite
{KearnsExact,DBLPconfnipsOrtizK02,DBLPconftarkGottlobGS03,ElkindGoldbergGoldbergEC06,DasPapEC06}}.

In the other direction, the notion of a PNE provides a \textit{new genre
of constraint satisfaction problems}; notably one in which, for any
assignment of strategies (values) to the neighborhood of a player
(variable), there is always a strategy (value) for that player which
makes the constraint (\ref{eq: equilibrium constraint}) corresponding
to that player satisfied (i.e., being in best response). The reason why
it might be hard to satisfy simultaneously the constraints
corresponding to all players is the {long-range correlations} that may
arise between players. Indeed, deciding whether a PNE exists is
NP-hard, even for very sparse graphical games \cite
{DBLPconftarkGottlobGS03}.

Viewed as a constraint satisfaction problem, the problem of the
existence of PNE poses interesting challenges. First, for natural
random ensembles over payoff tables such as the one adopted in this
paper (see Definition~\ref{def:model}), the \textit{expected number} of
PNE is $1$ for \textit{any graph} [this is shown for our model in the main
body of the paper; see (\ref{eqexpetation})].
On the contrary, for typical constraint satisfaction problems, the
expected number of solutions is exponential in the size of the graph
with different exponents corresponding to different density parameters.
Second, unlike typical constraint satisfaction problems studied before,
the existence of PNE is a priori not a monotone property of the connectivity.
It is surprising that given these novel features of the problem it is
possible to obtain a result establishing a double phase transition on
the existence of PNE as described below.


\subsection{Our model}

We define the notion of a graphical game and proceed to describe the
ensemble of random graphical games studied in this paper.
%
%
%
\begin{definition}[(Graphical game)] \label{def:graphical game}
Given a graph $G=(V,E)$, we define the \textit{neighborhood of node $v
\in
V$} to be the set $\N(v)= \{ v' | (v,v') \in E\}$. If $S_v$ is a set
associated with vertex $v$, for all $v \in V$, we denote by $S_{\N(v)}
:= \bigtimes_{v' \in\N(v)} S_{v'}$ the Cartesian product of the sets
associated with the nodes in $v$'s neighborhood.

A \textit{graphical game} on $G$ is a collection $(S_v, u_v)_{v \in V}$,
where $S_v$ is the \textit{strategy set of node $v$} and $u_v \dvtx S_v
\times
S_{\N(v)} \rightarrow\mathbb{R}$ the \textit{utility} (\textit{or payoff})
\textit{function} (\textit{or table}) \textit{of player $v$}. We also
define \textit{the best
response function} (\textit{or table}) \textit{of player $v$} to be the
function $\BR_v
\dvtx S_v \times S_{\N(v)} \rightarrow\{0,1\}$ such that
\[
\BR_v\bigl(\sigma_v, \sigma_{\N(v)}\bigr) = 1 \quad\Leftrightarrow\quad
\sigma_v \in
\arg_x \max\bigl\{ u_v\bigl(x, \sigma_{\N(v)}\bigr) \bigr\},
\]
for all $\sigma_v \in S_v$ and $\sigma_{\N(v)} \in S_{\N(v)}$.
\end{definition}
\begin{definition}[(Random graphical games on a fixed graph)] \label{def:model}
Given a graph $G=(V,E)$ and an atomless distribution $\F$
over $\mathbb{R}$, the probability distribution $\D_{G, \F}$ over
graphical games $(S_v, u_v)_{v \in V}$ on $G$ is defined as follows:
\begin{itemize}
\item$S_v = \{0,1\}$, for all $v\in V$;
\item the payoff values $\{u_v (\sigma_v, \sigma_{\N(v)}) \}_{v\in
V, \sigma_v \in S_v, \sigma_{\N(v)} \in S_{\N(v)}}$ are mutually
independent and identically distributed according to $\F$.
\end{itemize}
%
\end{definition}
%
%
\begin{remark}[(Invariance under payoff distributions)]\label
{remark:best response is enough}
It is easy to see that the existence of a PNE is only determined by the
\textit{best response tables} of the game; see condition (\ref{eq:
equilibrium constraint}). In particular, given that the distributions
considered in this paper are atomless, we can study PNE under $\D_{G,
\F}$, for any atomless $\F$, by restricting our attention (up to
probability $0$ events) to the measure $\D_G$ over best response
tables, defined as follows:
\begin{itemize}
\item$\{\BR_v (0, \sigma_{\N(v)}) \}_{v\in V, \sigma_{\N(v)} \in
S_{\N(v)}}$ are mutually independent and uniform in $\{0,1\}$;
\item$\BR_v (1, \sigma_{\N(v)}) =1-\BR_v (0, \sigma_{\N(v)})$,
for all $ \sigma_{\N(v)} \in S_{\N(v)}$.
\end{itemize}
We will sometimes refer to a graphical game defined in terms of its
best response tables as an \textit{underspecified graphical game}. Other
times, we will overload our terminology and just call it a graphical
game. We use $\PrGG[\cdot]$ and $\EGG[\cdot]$ to denote probabilities of
events and expectations, respectively, under the measure $\D_{G}$.
\end{remark}
\begin{remark}[(Invariance under payoff distributions II)]
\label{remark:different distributions for different players}
Given our observation in Remark~\ref{remark:best response is enough},
it follows that it is not important to use a common distribution $\F$
for sampling the payoffs of all the players of the game. All our
results in this work are true if different players have different
distributions as long as these distributions are atomless and all
payoffs values are sampled independently.
\end{remark}

\subsubsection*{Extending the model to random graphs}

One of the goals of this paper is to investigate what average degree is
required in a graph for a graphical game played on this graph to have a
PNE. To study this question, it is natural to consider families of
graphs with different densities and relate the probability of PNE
existence with the density of the graph. We consider graphical games on
graphs drawn from the Erd\H{o}s--R\'enyi, $G(n,p)$, model, with
varying values of the edge probability $p$. The ensemble of graphical
games we consider is formally the following.
\begin{definition}\label{def: random graphical games on random graphs}
Given $n \in\mathbb{N}$, $p \in[0,1]$ and an atomless distribution
$\F$ over~$\mathbb{R}$, we define a measure $\D_{(n,p, \F)}$ over
graphical games. A~graphical game is drawn from $\D_{(n,p, \F)}$ as follows:
\begin{itemize}
\item a graph $G$ is drawn from $G(n,p)$;
\item a random graphical game is drawn from $\D_{G, \F}$.
\end{itemize}
\end{definition}
\begin{remark}[(Invariance under payoff distributions III)]
\label{rem:random graphical game on random graphs}
Given our discussion in Remark~\ref{remark:best response is enough},
it follows that in order to study PNE in the random ensemble of
Definition~\ref{def: random graphical games on random graphs}, it is
sufficient to consider a measure that fixes only the best response
tables of the players in the sampled games.

For a given $n \in\mathbb{N}$ and $p \in[0,1]$, we define the
measure $\D_{(n,p)}$ over underspecified graphical games. An
underspecified graphical game is drawn from $\D_{(n,p)}$ as follows:
\begin{itemize}
\item a graph $G$ is drawn from $G(n,p)$;
\item a random underspecified graphical game is drawn from $\D_{G}$.
\end{itemize}
We use $\PrGT[\cdot]$ to denote probabilities of events under the measure
$\D_{(n,p)}$ and $\PrG[\cdot]$ for probabilities of events measurable
under $G(n,p)$.
\end{remark}

In the model defined in Definition~\ref{def: random graphical games on
random graphs} and Remark~\ref{rem:random graphical game on random
graphs}, there are two sources of randomness: the selection of the
graph, determining what players interact with each other, and the
selection of the payoff tables.
Note that in the two-stage process that samples a graphical game from
our distribution, the payoff tables can only be realized once the graph
is fixed. This justifies the subscript $G$ in the measure $\D_G$
defined above.

\subsection{Discussion}

\subsubsection*{Nonmonotonicity} Observe that the existence of a PNE is a
\textit{nonmonotone property} of $p$: any graphical game on the empty
graph has a PNE for trivial reasons; on the\vspace*{1pt} complete graph
a random
graphical game has a PNE with asymptotic probability $1-\frac{1}{e}$
(see~\cite{Dresher70,Rinott00}); but our results indicate that, when
$p$ is in some intermediate regime, a~PNE does not exist with
probability approaching $1$ as $n \rightarrow+\infty$ (see
Theorem~\ref{thm:mediumconn}).

The nonmonotonicity in the average degree of the existence of a PNE
makes the relation between PNE and connectivity nonobvious.
Surprisingly, we show (Theorem~\ref{th: high connectivity}) that the
convergence to a Poisson distribution of the distribution of the number
of PNE in complete graphs~\cite{Powers90,Stanford95} extends to much
sparser graphs, as long as the average degree is at least
\textit{logarithmic} in the number of players. If the sparsity increases
further, we show (Theorem~\ref{thm:mediumconn}) that a PNE does not
exist with high probability, while if the graph is essentially empty,
PNE exist with probability~$1$ (Theorem~\ref{thm:lowconn}). Our
results establish a \textit{double phase transition} consistent with the
nonmonotonicity of the model.

\subsubsection*{Methodological challenges} Our study here is an
instance of
studying the satisfiability of constraint satisfaction problems (CSPs).
The generic question is to investigate the effect of the structure of
the constraint graph on the satisfiability of the problems defined on
that graph, as well as their computational complexity.
In the context of CNF formulas (corresponding to the satisfiability
problem) the graph property most commonly studied in the literature is
the density of the hypergraph that contains an edge for each clause of
the formula (see, e.g.,~\cite{Friedgut99}). 
In other settings, different structural properties of the constraint
graph are relevant, for example, measures of cyclicity of the
graph \cite
{YannakakisAcyclicDBschemes,GottlobHypertreeDecompositions}. In our
case, studying the average degree reveals an interesting, nonmonotonic
behavior of the model, as described above.

In a typical CSP, to show that a solution does not exist one either
uses the first moment method to exhibit that the expected number of
solutions is tiny~\cite{AlonSpencer00}, or finds a witness of
unsatisfiability that exists with high probability.
To show that a satisfying assignment \textit{does} exist it is quite
common to use the second moment method or its refinements, which have
provided some of
the best bounds for satisfiability to date~\cite{achlioptasperes}.
In our model the expected number of satisfying assignments turns out to be
$1$ for any graph [see (\ref{eqexpetation}) below]. This
suggests that the analysis of the problem should be harder, since in
particular we cannot use\vadjust{\goodbreak} the first moment method to establish the
nonexistence of a PNE. Our proof of the nonexistence of PNE
(Theorems~\ref{thm:mediumconn} and~\ref{th:easy matching pennies})
uses succinct nonexistence witnesses that appear with high probability
in sufficiently sparse graphs. These witnesses are specific subgame
structures that do not possess a PNE with high probability. To
establish the existence of a PNE for sufficiently large densities
(Theorems~\ref{th: high connectivity} and~\ref{thm:expander}) we use
the second moment method.
Further, we use Stein's~\cite{Arratia89} method to establish that the
distribution of the number of PNE converges asymptotically to a
$\operatorname{Poisson}(1)$ distribution in this case.

\subsection{Outline of main results} \label{sec: main results}

We describe first our results for random graphs (for the measure $\D
_{(n,p)}$ defined in Remark~\ref{rem:random graphical game on random
graphs}), and proceed with our results for deterministic graphs (for
the measure $\D_G$ defined in Remark~\ref{remark:best response is enough}).

\subsubsection*{PNE on random graphs}

We study how the connectivity probability $p$ influences the existence
of PNE for games sampled from $\D_{(n,p)}$. The transition is
described by the following theorems applying to different levels of
graph connectivity. Before stating the theorems, we introduce some notation.
\begin{remark}[(Order notation)]
Let $f(x)$ and $g(x)$ be two functions
defined on some subset of the real numbers. One writes $f(x)=O(g(x))$
if and only if, for sufficiently large values of $x$, $f(x)$ is at most
a constant times $g(x)$ in absolute value. That is, $f(x) = O(g(x))$ if
and only if there exists a positive real number $M$ and a real number
$x_0$ such that
\[
|f(x)| \le M |g(x)|\qquad\mbox{for all }x>x_0.
\]
Similarly, we write $f(x) = \Omega(g(x))$ if and only if there exists
a positive real number $M$ and a real number $x_0$ such that
\[
|f(x)| \ge M |g(x)|\qquad\mbox{for all }x>x_0.
\]

We casually use the \textit{order notation} $O(\cdot)$ and $\Omega(\cdot
)$ throughout the paper. Whenever we use $O(f(n))$ or $\Omega(f(n))$
in some bound, there exists a constant $c >0 $ such that the bound
holds true for sufficiently large $n$ if we replace the $O(f(n))$ or
$\Omega(f(n))$ in the bound by $c \cdot f(n)$.
\end{remark}
\begin{remark}[(Order notation continued)]
If $g(n)$ is a function of $n \in\mathbb{N}$, then we denote by
$\omega(g(n))$ any function $f(n)$ such that ${f(n) / g(n)}
\rightarrow+\infty$, as $n \rightarrow+\infty$; similarly, we
denote by $o(g(n))$ any function $f(n)$ such that ${f(n) / g(n)}
\rightarrow0$, as $n \rightarrow+\infty$. Finally, for two functions
$f(n)$ and $g(n)$, we write $f(n) \gg g(n)$ whenever $f(n) = \omega(g(n))$.
\end{remark}
\begin{theorem}[(High connectivity)] \label{th: high connectivity}
Let $Z$ denote the number of PNE in
a graphical game sampled from $\D_{(n,p)}$,
where $p=\frac{(2+\varepsilon) \log_e (n)}{n}$, $\eps= \eps(n) > 0$.
For an arbitrary constant $c>0$ we assume that $\eps(n) >c$ and (in
order for $p \le1$) $\eps(n) \le{n \over\log_e(n)}-2$.\vadjust{\goodbreak}


Under the above assumptions, for all finite $n$, with probability at
least $1 - 2 n^{-\eps/8}$ over the random graph sampled from $G(n,p)$,
it holds that
the total variation distance between $Z$ and a $\operatorname
{Poisson}(1)$ r.v.
$W$ is bounded by
%
%
\begin{equation}\label{eq: high connectivity total variation distance}
\| Z - W \| \leq O(n^{-\varepsilon/8}) + \exp(-\Omega(n)).
\end{equation}
In other words,
%
%
\begin{equation}
\PrG[ \| Z - W \| \leq O(n^{-\varepsilon/8}) + \exp(-\Omega(n))
]
\geq1 - 2 n^{-\eps/8}.
\end{equation}


In particular, the distribution of $Z$ converges in total variation
distance to a $\operatorname{Poisson}(1)$
distribution, as $n \rightarrow+\infty$.

[Note that the two terms on the right-hand side
of (\ref{eq: high connectivity total variation distance}) can be of
the same order when $\eps$ is of the order of $n/ \log_e (n)$.]
\end{theorem}
%
%
%
\begin{theorem}[(Medium connectivity)]
\label{thm:mediumconn}
For all $p=p(n) \leq1/n$, if a graphical game is sampled from $\D
_{(n,p)}$, the probability that a PNE exists is bounded by
\[
\exp(-\Omega(n^2 p)).
\]
For $p(n) = g(n)/n$, where $\log_e (n)/2>g(n) > 1$, the probability
that a PNE exists is bounded by
\[
\exp\bigl(-\Omega\bigl(e^{\log_e (n) -2 g(n)}\bigr)\bigr).
\]
In particular, the probability that a PNE exists goes to $0$ as $n
\rightarrow+\infty$ for all $p=p(n)$ satisfying
\[
\frac{1}{n^2} \ll p < \bigl(0.5 - \varepsilon'(n) \bigr) \frac{\log_e (n)}{n},
\]
where $\varepsilon'(n) = \omega( {1\over\log_e (n)})$.
\end{theorem}
\begin{theorem}[(Low connectivity)]
\label{thm:lowconn}
For every constant $c>0$, if a graphical game is sampled from $\D
_{(n,p)}$ with $p \le\frac{c}{n^2}$, the probability that a PNE exists
is at least
\[
\biggl(1-\frac{c}{n^2}\biggr)^{{n (n-1)}/{2}} \longrightarrow
e^{-{c/2}}.
\]
\end{theorem}

Note that our upper and lower bounds for $G(n,p)$ leave a small gap,
between $p \approx\frac{0.5 \log_e (n)}{n}$
and $p \approx\frac{2 \log_e (n)}{n}$. The behavior of the number of
PNE in this range of $p$ remains open. We establish the nonexistence
of PNE for medium connectivity graphs via a simple structure that
prevents PNE from arising, called the ``indifferent matching pennies
game'' (see Definition~\ref{def:indifferent matching pennies} below).
It is natural to ask whether our ``indifferent matching pennies''
witnesses are (similarly to isolated vertices in connectivity) the
smallest structures that prevent the existence of PNE and the last ones
to disappear.

%

\subsubsection*{General graphs}
We give conditions on the structure of a graph implying the (likely)
existence or nonexistence of a PNE in a random game played on that
graph. The existence
of a PNE is guaranteed by sufficient connectivity
of the underlying graph. The connectivity that we require is captured
by the notion of \textit{$(\alpha,\delta)$-expansion} given next.
\begin{definition}[{[$(\alpha,\delta)$-expansion]}]
A graph $G=(V,E)$ has \textit{$(\alpha,\delta)$-expan\-sion} iff
every set $V'$ such that $|V'| \le\lceil\delta|V| \rceil
$ has $|\mathcal{N}(V')|
\ge\min(|V|, \alpha|V'|)$ neighbors.
Here we let
\[
\mathcal{N}(V') = \{ w \in V \dvtx\exists u \in V' \mbox{ with }
(u,w) \in E\}.
\]
[Note in particular that $\mathcal{N}(V')$ may intersect $V'$.]
\end{definition}

We show the following result.
\begin{theorem}[(Strongly connected graphs)]\label{thm:expander}
Let $Z$ denote the number of PNE in a graphical game sampled from $\D
_G$, where $G$ is a graph on $n$ vertices
that has $(\alpha,\delta)$-expansion with
$\alpha=(1+\varepsilon)\log_2 (n)$, $\delta=\frac{1}{\alpha}$ and
\mbox{$\varepsilon>0$}. Then
the total variation distance between the distribution of $Z$ and the
distribution of a $\operatorname{Poisson}(1)$ r.v. $W$ is bounded by
%
%
\begin{equation}
\| Z - W\| \leq O(n^{-\varepsilon}) + O(2^{-n/2}).
\end{equation}
\end{theorem}

Next we provide a complementary condition for the
nonexistence of PNE. The condition will be given in terms of the
following structure.
\begin{definition}[($d$-bounded edge)]
An edge $e =(u,v) \in E$ of a graph $G(V,E)$ is called \textit
{d-bounded} if both $u$ and $v$ have degrees smaller or equal to~$d$.
\end{definition}

We bound the probability that a PNE exists in a game sampled from $\D
_G$ as a function of the number of $d$-bounded edges in $G$. For the
stronger version of our theorem, we also need the notion of a
\textit{maximal weighted independent edge-set} defined next.
\begin{definition}[(Maximal weighted independent edge-set)]
Given a graph $G(V,E)$, a~subset $\E\subseteq E$ of the edges is
called \textit{independent} if no pair of edges in $\E$ are adjacent.

If $w\dvtx E \rightarrow\mathbb{R}$ is a function assigning weights to
the edges of $G$, we extend $w$ to subsets of edges by assigning to
each $\E\subseteq E$ the weight $w_{\E} = \sum_{e \in\E} w(e)$.
Then we call a subset $\E\subseteq E$ of edges \textit{a maximal weighted
independent edge-set} if $\E$ is an independent edge-set with maximal
weight among independent edge-sets.
\end{definition}
\begin{theorem} \label{th:easy matching pennies}
A random game sampled from $\D_G$, where $G$ is a graph with at least
$m$ vertex disjoint $d$-bounded edges, has no PNE with probability\vadjust{\goodbreak} at least
%
%
\begin{equation} \label{eq:edge_disjoint1}
1-\exp\bigl(-m \bigl(\tfrac{1}{8}\bigr)^{2^{2d-2}}\bigr).
\end{equation}
In particular, if $G$ has at least $m$ edges that are $d$-bounded,
then a game sampled from $\D_G$ has no PNE with probability at least
%
%
\begin{equation} \label{eq:edge_disjoint2}
1-\exp\biggl(-\frac{m}{2d} \biggl(\frac{1}{8}
\biggr)^{2^{2d-2}}\biggr).
\end{equation}
Moreover, there exists an algorithm of complexity
$O(n^2 + m 2^{d+2})$
%
for proving that a PNE does not exist, which has success probability
given by (\ref{eq:edge_disjoint1}) and~(\ref{eq:edge_disjoint2}),
respectively.

More generally, let us assign to every edge $(u,v) \in E$ the weight
\[
w_{(u,v)}:= - \log_e\bigl(1-p_{(u,v)}\bigr),
\]
for $p_{(u,v)} = 8^{-2^{d_u+d_v-2}}$, where $d_u$ and $d_v$ are,
respectively, the degrees of $u$ and~$v$.
Given these weights, suppose that $\E$ is a maximal weighted
independent edge-set with
value $w_{\E}$. Then the probability that there exists no PNE is at least
\[
1-\exp(-w_{\E}).
\]
\end{theorem}

An easy consequence of this result is that many sparse graphs, such as
the line and the grid, do not have a PNE with probability tending to
$1$ as the number of players increases.

The proof of Theorem~\ref{th:easy matching pennies} is based on a
small witness for the nonexistence of PNE, called the \textit{indifferent
matching pennies game}. As the name implies, this game is inspired by
the simple \textit{matching pennies game}. Both games are described next.
\begin{definition}[(The matching pennies game)] \label{def:matching
pennies game}
We say that two players $a$ and $b$ play the \textit{matching pennies
game} if their payoff matrices are the following, up to permuting the
players' names.
\[
\begin{tabular}{@{}cccc@{}}
\multicolumn{4}{@{}c@{}}{Payoff table of player $a$:}\\
& \vline & $b$ plays 0 & $b$ plays 1 \\ \hline
$a$ plays 0 & \vline& 1 & 0 \\
$a$ plays 1 & \vline& 0 & 1
\end{tabular}
\hspace*{20pt}
\begin{tabular}{@{}cccc@{}}
\multicolumn{4}{@{}c@{}}{Payoff table of player $b$:}\\
& \vline& $b$ plays 0 & $b$ plays 1 \\ \hline
$a$ plays 0 & \vline& 0 & 1 \\
$a$ plays 1 & \vline& 1 & 0
\end{tabular}
\]
\end{definition}
\begin{definition}[(The indifferent matching pennies game)]
\label{def:indifferent matching pennies}
We say that two players $a$ and $b$ that are adjacent to each other in
a graphical game play the \textit{indifferent matching pennies game} if,
for all strategy profiles $\sigma_{\N(a) \cup{\N(b)}\setminus\{
a,b\}}$ in the neighborhood of $a$ and $b$, the players $a$ and $b$
play a matching pennies game against each other.

In other words, for all fixed $\sigma:=\sigma_{\N(a) \cup{\N
(b)}\setminus\{a,b\}}$, the payoff tables of $a$ and $b$ projected on
$\sigma_{\N(a)\setminus\{b\}}$ and $\sigma_{\N(b)\setminus\{a\}
}$, respectively,\vadjust{\goodbreak} are the following, up to permuting the players' names.
\begin{eqnarray*}
\begin{tabular}{@{}lccc@{}}
\multicolumn{4}{@{}c@{}}{Payoffs to player $a$:}\\[4pt]
& \vline&
\begin{minipage}{4 cm}
\centering{\mbox{$b$ plays 0, other neighbors} \mbox{play $\sigma_{\N
(a)\setminus
\{b\}}$}}
\end{minipage}
&
\begin{minipage}{4 cm}\centering{\mbox{$b$ plays 1, other neighbors}
\mbox{play $\sigma_{\N(a)\setminus\{b\}}$}}
\end{minipage}
\\ \hline
$a$ plays 0 & \vline& 1 & 0 \\
$a$ plays 1 & \vline& 0 & 1
\end{tabular}
\\
\begin{tabular}{@{}lccc@{}}
\multicolumn{4}{@{}c@{}}{Payoffs to player $b$:}\\[4pt]
& \vline&
\begin{minipage}{4.1 cm} \centering{\mbox{$a$ plays 0, other neighbors}
\mbox{play $\sigma_{\N(b)\setminus
\{a\}}$}}
\end{minipage}
&
\begin{minipage}{4.1 cm} \centering{\mbox{$a$ plays 1, other neighbors}
\mbox{play $\sigma_{\N(b)\setminus
\{a\}}$}}
\end{minipage}
\\ \hline
$b$ plays 0 & \vline& 0 & 1 \\
$b$ plays 1 & \vline& 1 & 0
\end{tabular}
\end{eqnarray*}
\end{definition}

Observe that if a graphical game contains an edge $(u,v)$ so that
players $u$ and $v$ play the indifferent matching pennies game, then
the game has no PNE.
In particular, the indifferent matching pennies game provides a \textit
{small witness} for the nonexistence of a PNE, which is a
coNP-complete problem for bounded degree graphical games \cite
{DBLPconftarkGottlobGS03}.
Our analysis implies that, with high probability over bounded degree
graphical games, there are short proofs for the nonexistence of PNE
which can be found efficiently.
A related analysis and randomized algorithm was introduced for mixed
Nash equilibria in two-player games by B\'ar\'any, Vempala and Vetta
\cite{Barany05}.




%



\subsection{Related work} The number of PNE in random games with
i.i.d. payoffs has been extensively studied in the literature prior to
our work: Goldberg, Goldman and Newman~\cite{Goldberg68} characterize the
probability that a two-player random game with i.i.d. payoff tables has
a PNE, as the number of strategies tends to infinity.
Dresher~\cite{Dresher70} and Papavassilopoulos~\cite{Papav95}
generalize this result to $n$-player random games on the complete graph.
Powers~\cite{Powers90} and Stanford~\cite{Stanford95} generalize the
result further, showing that the distribution of the number of PNE
approaches a $\operatorname{Poisson}(1)$ distribution
as the number of strategies increases. Finally, Rinott and Scarsini
\cite
{Rinott00} investigate the asymptotic distribution of PNE for a more
general ensemble of random games on the complete graph where there are
positive or negative dependencies among the players' payoffs.

Our work generalizes the above results for i.i.d. payoffs beyond the
complete graph to random graphical games on random graphs and several
families of deterministic graphs. Parallel to our work, Dilkina, Gomes
and Sabharwal
\cite{dilkinaaaai07} studied the existence of PNE in certain
families of deterministic graphs, and Hart, Rinott and Weiss~\cite{HRW08}
obtained results for evolutionarily stable strategies in random games.
These results are related but not directly comparable to our results.


\section{Random graphs}
\subsection{High connectivity}
In this section we study the number of PNE in graphical games sampled
from $\D_{(n,p)}$. We show that, when the average degree is $pn =
(2+\varepsilon(n)) \log_e (n)$, where $\varepsilon(n) >c$ and $c>0$ is any
fixed constant, the distribution of the number of PNE converges to a
$\operatorname{Poisson}(1)$ random variable, as $n$ goes to infinity.
This implies in particular that a
PNE exists with probability converging to
$1-\frac{1}{e}$ as the size of the network increases.

As in the study of the complete graph in~\cite{Rinott00}, we use the
following result of Arratia, Goldstein and Gordon~\cite{Arratia89},
established using
Stein's method.
For two random variables $Z,Z'$ supported on $0,1,\ldots$ we define
their \textit{total variation distance} $\| Z - Z' \|$
as
\[
\| Z - Z' \| := \sum_{i=0}^{\infty} |Z(i) - Z'(i)|.
\]

%
\begin{lemma}[\cite{Arratia89}] \label{lem:Arratia}
Consider arbitrary Bernoulli random variables $X_i$, $i=0,\ldots,N$.
For each $i$, define some \textit{neighborhood of dependence} $B_i$ of $X_i$
such that $B_i$ satisfies that
$(X_j \dvtx j \in B^c_i)$ are \textit{independent} of $X_i$.
Let
%
%
\begin{equation}
Z = \sum_{i=0}^{N} X_i,\qquad \lambda= \mathbb{E}[Z],
\end{equation}
and assume that $\lambda>0$. Also, let
\[
b_1= \sum_{i=0}^{N} \sum_{j \in B_i} \Prob[X_i=1]\Prob[X_j=1]
\]
and
\[
b_2= \sum_{i=0}^{N} \sum_{j \in B_i \setminus\{i\}}
\Prob[X_i=1,X_j=1].
\]
%
Then the total variation distance between the distribution of $Z$
and a $\mathrm{Poisson}$ random variable $W_\lambda$ with mean $\lambda$
is bounded by
%
%
\begin{equation}
\|Z - W_\lambda\| \leq2(b_1+b_2).
\end{equation}
\end{lemma}
%
%
\begin{pf*}{Proof of Theorem~\ref{th: high connectivity}}
For ease of notation, we identify the players of the
graphical game with the indices $1, 2,\ldots,n$. We also identify pure
strategy profiles with the integers in $\{0,\ldots, 2^n-1\}$, mapping
each integer to a strategy profile. The mapping is defined so that, if the
binary expansion of $i$ is $i(1) \cdots i(n)$, player $k$ plays $i(k)$.

Next, to each strategy profile $i \in\{0, \ldots, N\}$, where
$N=2^n-1$, we assign an indicator random variable $X_{i}$ which is $1$
if the strategy profile $i$ is a PNE.\vadjust{\goodbreak} Then the counting random variable
%
%
\begin{equation}
Z=\sum_{i=0}^N X_{i}
\end{equation}
corresponds to the number of PNE. Hence the existence of a PNE is
equivalent to the random variable $Z$ being positive.

Let us condition on a realization of the graph $G$ of the graphical
game, but not its best response tables. For a given strategy profile
$i$, each player is in best response with probability $1/2$ over the
selection of her best response table\setcounter{footnote}{3}\footnote
{This follows directly
from our model (Remark~\ref{rem:random graphical game on random
graphs}), following our assumption of atomless payoff distributions
(Definition~\ref{def: random graphical games on random graphs}).};
therefore $\stexp_G[X_i] = 2^{-n}$, for all $i$, where we recall that
$\stexp_G$ denotes expectation under the measure $\D_G$.
Hence, conditioning on $G$ the expected number of PNE is
%
%
\begin{equation}
\label{eqexpetation}
\stexp_G[Z] = 1.
\end{equation}
Since this holds for any realization of the graph $G$ it follows that
$\stexp[Z] = 1$.


%

In Lemma~\ref{lem:neighborhood of dependence} that follows, we
characterize the neighborhood of dependence $B_i$
of the variable $X_i$ in order to be able to apply Lemma \ref
{lem:Arratia} on the collection of variables $X_0, \ldots, X_N$. Note
that this neighborhood depends on the graph realization, but is
independent of the realization of the payoff tables.
\begin{lemma}\label{lem:neighborhood of dependence} For a fixed graph
$G$, we can choose the neighborhoods of dependence for the random
variables $X_0, \ldots, X_N$ as follows:
\[
B_0 = \{ j \dvtx\exists k \mbox{ such that } \forall k' \mbox{ with }
(k,k') \in E(G)
\mbox{ it holds that } j(k') = 0\}
\]
and
\[
B_i = i \oplus B_0 = \{ i \oplus j \dvtx j \in B_0\},
\]
where $i \oplus j = (i(1) \oplus j(1),\ldots,i(n) \oplus j(n))$ and
$\oplus$ is the exclusive or operation.
\end{lemma}
\begin{remark}
Intuitively, when the graph $G$ is realized, the neighborhood of
dependence of the strategy profile $0$ (variable $X_0$) contains all
strategy profiles $j$ (variables $X_j$) assigning $0$ to all the
neighbors of at least one player $k$. If such a player $k$ exists, then
whether $0$ or $j(k)$ is a best response to the all-$0$ neighborhood
are dependent random variables (over the selection of the best response
table of player $k$). The definition of $B_i$ in terms of $B_0$ is
justified by the symmetry of our model.
\end{remark}
\begin{pf*}{Proof of Lemma~\ref{lem:neighborhood of dependence}}
By symmetry, it is enough to show that $X_0$ is independent of $\{X_i\}
_{i \notin B_0}$. Fix some $i \notin B_0$. Observe that in $i$, each
player $k$ of the game has at least one neighbor $k'$ playing strategy
$1$. By the definition of measure $\D_G$,\vadjust{\goodbreak} it follows that whether
strategy $0$ is a best response for player $k$ in strategy profile $0$
is independent of whether strategy $i(k)$ is a best response for player
$k$ in strategy profile $i$, since these events depend on different
strategy profiles of the neighbors of $k$.
\end{pf*}

Now, for a fixed graph $G$, the functions $b_1(G)$ and $b_2(G)$
(corresponding to $b_1$ and $b_2$ in Lemma~\ref{lem:Arratia}) are well defined.
We proceed to bound the expectation of these functions over the
sampling of the graph $G$:
%
%
\begin{eqnarray}
\stexp_{\G}[ b_1 (G)] &=& \stexp_{\G}\Biggl[ \sum_{i=0}^{N} \sum
_{j \in
B_i} \PrGG[X_i=1]\PrGG[X_j=1]\Biggr] \nonumber\\
&=& \stexp_{\G}\Biggl[
\frac{1}{(N+1)^2} \sum_{i=0}^{N} |B_i|\Biggr] \\
&=&
\frac{\stexp_{\G}[|B_0|]}{N+1};
\nonumber\\
\label{eq: expectation of b2}
\stexp_{\G}[ b_2(G)]&=& \stexp_{\G}\Biggl[ \sum_{i=0}^{N} \sum_{j
\in
B_i \setminus\{i\}}\PrGG[X_i=1,X_j=1]\Biggr] \nonumber\\[-8pt]\\[-8pt]
&=& (N+1) \sum
_{j \neq0}
\stexp_{\G}\bigl[ \PrGG[X_0=1,X_j=1] \ones[j \in B_0]\bigr].\nonumber
\end{eqnarray}
In the last line of both derivations we made use of the symmetry of the
model. Invoking symmetry again, we observe that the expectation
\[
\stexp_{\G}\bigl[ \PrGG[X_0=1,X_j=1] \ones[j \in B_0]\bigr]
\]
in (\ref{eq: expectation of b2}) depends only on the number of $1$'s
in the strategy profile $j$, denoted $s$ below.
Let us write $Y_s$ for the indicator that the strategy profile $j_s$,
where the first $s$ players play $1$, and all the other players play
$0$, is a PNE. Also, write $I_s$ for the indicator that this strategy
is in $B_0$ (note that $I_s$ is a function of the graph only). Using
this notation, we obtain
%
%
\begin{equation} \label{eq:b1}
\stexp_{\G} [b_2(G)]= 2^n\sum_{s=1}^{n} \pmatrix{n \cr s}
\stexp_{\G}\bigl[ I_s \PrGG[ Y_0=1, Y_s=1] \bigr]
\end{equation}
and
\begin{equation}
\stexp_{\G} [b_1(G)]= 2^{-n} \sum_{s=0}^{n} \pmatrix{n \cr s}
\stexp_{\G}[I_s].
\end{equation}

\begin{lemma} $\stexp_{\G} [b_1(G)]$ and $\stexp_{\G} [b_2(G)]$ are
bounded as follows:
\begin{eqnarray*}
\stexp_{\G} [b_1(G)] &\leq&
R(n,p) :=\sum_{s=0}^n \pmatrix{n \cr s} 2^{-n} \min\bigl(1, n
(1-p)^{s-1}\bigr),\\
\stexp_{\G} [b_2(G)] &\leq&
S(n,p) :=\sum_{s=1}^n \pmatrix{n \cr s} 2^{-n}
\bigl[\bigl(1+(1-p)^s\bigr)^{n-s} -
\bigl(1-(1-p)^s\bigr)^{n-s} \bigr].
\end{eqnarray*}
\end{lemma}
\begin{pf}
We begin with the study of $\stexp_{\G} [b_1(G)]$. Clearly, it
suffices to bound $\stexp[I_s]$ by $n (1-p)^{s-1}$, for $s>0$. For the
strategy profile $j_s$ to belong in $B_0$ it must be that there is at
least one player who is not connected to any player in the set $S:=\{
1,2,\ldots,s\}$.
The probability that a specific player $k$ is not connected to any
player in $S$ is either $(1-p)^s$ or $(1-p)^{s-1}$, depending on
whether $k \in S$; so it is always at most $(1-p)^{s-1}$. By a union
bound it follows that the probability there is at least one player not
connected to $S$ is at most $n (1-p)^{s-1}$.

We now analyze $\stexp_{\G}[ I_s \PrGG[ Y_0=1, Y_s=1] ]$.
Recall from the previous paragraph that $I_s = 1$ only when there
exists a player $k$ who is not connected to any player in~$S$. If there
exists such a player $k$ with the extra property that $k \in S$, then
$\PrGG[Y_0=1 , Y_s=1]=0$, since it cannot be that both $0$ and $1$ are
best responses for player $k$ when all her neighbors play $0$.

Therefore the only contribution to $\stexp_{\G}[ I_s \PrGG[ Y_0=1,
Y_s=1] ]$
is from the event every player in $S$ is connected to at least one
other player in $S$. 
Conditioning on this event, in order for $I_s = 1$ it must be that
at least one of the players in $S^c:= V\setminus S$ is not adjacent to
any player in $S$.

Let us define $p_s:=\PrG[ \not\exists\mbox{ isolated node in the
subgraph induced by $S$}]$, and let
$t$ denote the number of players in $S^c$, which are not connected to
any player in~$S$.
Since every player outside $S$ is nonadjacent to any player in $S$
with probability $(1-p)^s$, the probability that exactly $t$ players
are not adjacent to $S$ is
\[
\pmatrix{n-s \cr t} [(1-p)^s]^t \bigl(1-(1-p)^s\bigr)^{n-s-t}.
\]
Moreover, conditioning on the event that exactly
$t$ players in $S^c$ are not adjacent to any player in $S$, 
we have that the probability that $Y_0 = 1$ and $Y_s = 1$ is
\[
\frac{1}{2^t} \frac{1}{2^{n-t}} \frac{1}{2^{n-t}}.
\]
Putting these together we obtain
\begin{eqnarray*}
&&\stexp_{\G}\bigl[ I_s \PrGG[ Y_0=1, Y_s=1] \bigr] \\
&&\qquad = p_s \sum_{t=1}^{n-s} \pmatrix{n-s \cr t} [(1-p)^s]^t
\bigl(1-(1-p)^s\bigr)^{n-s-t} \frac{1}{2^t} \frac{1}{4^{n-t}} \\
&&\qquad = \frac{p_s}{4^n} \bigl( \bigl( 2 (1-p)^s +
\bigl(1-(1-p)^s\bigr) \bigr)^{n-s} -
\bigl(1 - (1-p)^s\bigr)^{n-s} \bigr) \\
&&\qquad
= \frac{p_s}{4^n} \bigl( \bigl(1 + (1-p)^s\bigr)^{n-s} - \bigl(1
- (1-p)^s\bigr)^{n-s} \bigr);
\end{eqnarray*}
therefore
\begin{eqnarray*}
\stexp_{\G} [b_2(G)]&=& \sum_{s=1}^{n} 2^{-n} \pmatrix{n \cr s} p_s
\bigl[ {
\bigl(1+(1-p)^s\bigr)^{n-s}- \bigl(1-(1-p)^s\bigr)^{n-s}
} \bigr] \\
&\leq& S(n,p).
\end{eqnarray*}
\upqed\end{pf}

In the \hyperref[app]{Appendix} we show  the following.
\begin{lemma} \label{lem:bounding expectations of b1 and b2}
\[
S(n,p) \leq O(n^{-\varepsilon/4}) + \exp(-\Omega(n))
\]
and
\[
R(n,p) \leq O(n^{-\varepsilon/4}) + \exp(-\Omega(n)).
\]
\end{lemma}

Given the above bounds on $\stexp_{\G} [b_1(G)]$ and $\stexp_{\G}
[b_2(G)]$, Markov's inequality implies that with probability at least
$1-n^{-\varepsilon/8} -2^{-n}$ over the selection of the graph $G$ from $G(n,p)$
we have
%
%
\begin{equation}\label{eq: conditionnnn}
\max(b_1(G),b_2(G)) \leq O(n^{-\varepsilon/8}) + \exp(-\Omega(n)).
\end{equation}
Let us condition on the event that condition (\ref{eq: conditionnnn})
holds. Under this event, Lem\-ma~\ref{lem:Arratia} implies that
\[
\| Z - W \| \leq2\bigl(b_1(G)+b_2(G)\bigr) \leq O(n^{-\eps/8}) + \exp(-\Omega(n))
\]
as needed. Noting that $1-n^{-\varepsilon/8} -2^{-n} \ge
1-2n^{-\varepsilon
/8}$, we obtain
\begin{equation}
\PrG[ \| Z - W \| \leq O(n^{-\varepsilon/8}) + \exp(-\Omega(n)) ] \geq1
- 2 n^{-\eps/8}.
\end{equation}

Using the pessimistic upper bound of $2$ on the total variation
distance when condition (\ref{eq: conditionnnn}) fails, we obtain
\[
\| Z - W \| \le O(n^{-\varepsilon/8}) + \exp(-\Omega(n)).
\]
Taking the limit of the above bound as $n \rightarrow+\infty$ we
obtain our asymptotic result. This concludes the proof of Theorem
\ref{th: high connectivity}.
\end{pf*}


\subsection{Medium connectivity}


\mbox{}

\begin{pf*}{Proof of Theorem~\ref{thm:mediumconn}}
Recall the \textit{matching pennies game} from Definition \ref
{def:matching pennies game}. 
It is not hard to see that this game does not have a PNE.
Hence, if a graphical game contains two players who are connected to
each other, are isolated from all the other players and play matching
pennies against each other, then the graphical game will have no PNE.
The existence of such a witness for the nonexistence of PNE is
precisely what we use to establish our result. In particular, we show
that with high probability a random game sampled from $\D_{(n,p)}$
will contain an isolated edge between two players playing a matching
pennies game.\vadjust{\goodbreak}

We use the following exposure argument. Label the vertices of the graph
with the integers in $[n]:=\{1,\ldots,n\}$. Set $\Gamma_1 = [n]$ and
perform the following operations, which iteratively define the sets of
vertices $\Gamma_i$, $i \ge2$. If\vspace*{1pt} $|\Gamma_i| \le n/2$, for some $i
\ge2$, stop the process and do not proceed to iteration $i$\footnote
{Throughout the process $\Gamma_i$ represents the set of vertices that
could be adjacent to an isolated edge, given the information available
to the process at the beginning of iteration~$i$.}:
\begin{itemize}
\item
Let $j$ be the minimal value such that $j \in\Gamma_i$.
\item
If $j$ is adjacent to more than one vertex or to none, let
$\Gamma_{i+1} = \Gamma_i \setminus(\{j\} \cup\N(j))$. Go to the
next iteration.
\item
Otherwise, let $j'$ be the unique neighbor of $j$. If $j'$ has a
neighbor $\neq j$, let $\Gamma_{i+1} = \Gamma_i \setminus(\{ j,j'\}
\cup\N(j'))$. Go to the next iteration.
\item
Otherwise check if the players $j$ and $j'$ play a matching pennies
game.\footnote{More precisely, check if the best response tables of
the players $j$ and $j'$ are the same with the best response tables of
the players $a$ and $b$ of the matching pennies game from
Definition~\ref{def:matching pennies game} (up to permutations of the
players' names).}
If this is the case, declare \textsc{No Nash}.
Let $\Gamma_{i+1} = \Gamma_i \setminus\{j, j'\}$. Go to the next iteration.
\end{itemize}

Observe that the number of vertices removed at some iteration of the
process can be upper bounded (formally, it is stochastically dominated) by
\[
2 + \operatorname{Bin}(n,p),
\]
where $\operatorname{Bin}(n,p)$ is a random variable distributed
according to
the binomial distribution with $n$ trials and success probability $p$.
This follows from the fact that the vertices removed at some iteration
of the process are either the examined vertex $j$ and $j$'s neighbors
[the number of those is stochastically dominated by a $\operatorname{Bin}(n,p)$
random variable], or---if $j$ has a single neighbor $j'$---the removed
vertices are $j$, $j'$ and the neighbors of $j'$ [the number of those
is also stochastically dominated by a $\operatorname{Bin}(n,p)$ random
variable]. Letting $m:= \lceil0.02 n / (np + 1) \rceil$, the
probability that the process runs for at most $m$ iterations is bounded by
\[
\Pr[ 2 m + \operatorname{Bin}(m n,p) \ge n/2 ] \le\exp(-\Omega(n)).
\]

Condition on the information known to the exposure process up until the
beginning of iteration $i$, and assume that $|\Gamma_i| > n/2$. Let
$j$ be the vertex with the smallest value in $\Gamma_i$. Now reveal
all the neighbors of $j$, and if $j$ has only one neighbor $j'$ reveal
also the neighbors of $j'$. The probability that $j$ is adjacent to a
node $j'$ who has no other neighbors is at least ${n \over2} p
(1-p)^{2n}=:p_{\mathrm{iso}}$; note that we made use of the condition
$|\Gamma_i| > n/2$ in this calculation. Conditioning on this event,
the probability (over the selection of the payoff tables) that $j$ and
$j'$ play a matching pennies game is ${1 \over8}=:p_{\mathrm{mp}}$. Hence,
the probability of outputting \textsc{No Nash} in iteration $i$ is at
least $\frac{1}{8} {1 \over2} n p (1-p)^{2n}=:p_{\mathrm{imp}}$.

%

The probability that the game has a PNE is upper bounded by the
probability that the process described above does not return \textsc{No
Nash}, at any point through its completion. To upper bound the latter
probability, let us imagine the following alternative process:
\begin{enumerate}
\item\textit{Stage} 1: Toss $n$ coins independently at random with head
probability $p_{\mathrm{iso}}$. Let $\I_1, \I_2,\ldots, \I_n \in\{0,1\}
$, where $1$ represents ``heads,'' and $0$ represents ``tails,'' be the
outcomes of these coin tosses.

\item\textit{Stage} 2: Toss $n$ coins independently at random with head
probability $p_{\mathrm{mp}}$. Let $\M_1, \M_2,\ldots, \M_n \in\{0,1\}
$, be the outcomes of these coin tosses.

\item\textit{Stage} 3: Run through the exposure process in the following
way. At each iteration $i$:
\begin{itemize}
\item conditioning on the information available to the exposure process
at the beginning of the iteration, compute the probability $p_j$ that
the vertex $j$ corresponding to the smallest number in $\Gamma_i$ is
adjacent to an isolated edge; given the discussion above it must be
that $p_j \ge p_{\mathrm{iso}}$;
\item if $\I_i =1$, then create an isolated edge connecting the player
$j$ to a random vertex $j' \in\Gamma_i \setminus\{j\}$, forbidding
all other edges from $j$ or $j'$ to any other player, and make the
players $j$ and $j'$ play a matching pennies game if $\M_i = 1$; if
they do output \textsc{No Nash};

\item if $\I_i =0$, then sample the neighborhood of $j$ from the
following modified model:
\begin{itemize}
\item with probability $p_j - p_{\mathrm{iso}} \over1 - p_{\mathrm{iso}}$,
create an isolated edge connecting the player $j$ to a random vertex
$j' \in\Gamma_i \setminus\{j\}$, forbidding all other edges from $j$
or $j'$ to any other player, and make the players $j$ and $j'$ play a
matching pennies game with probability $p_{\mathrm{mp}}$; if both of these
happen, output \textsc{No Nash};

\item with the remaining probability, sample the neighborhood of $j$
and the neighborhood of the potential unique neighbor $j'$ from
$G(n,p)$, conditioning on $j$ not being adjacent to an isolated edge.
\end{itemize}
\item Define $\Gamma_{i+1}$ from $\Gamma_i$ appropriately, and exit
the process if $|\Gamma_{i+1}| \le n/2$.
\end{itemize}
\end{enumerate}
It is clear that the process given above can be coupled with the
process defined earlier to exhibit the same behavior. But it is easier
to analyze. In particular, letting $\mathcal{S}:= \sum_{i=1}^m \I_i
\M_i $, the probability that a Nash equilibrium does not exist can be
lower bounded as follows:
\begin{eqnarray*}
\PrG[\not\exists\mbox{ a PNE}] &\ge&\Pr\left[ \mathcal{S} \ge
1 \wedge
\begin{array}{c}
\mbox{process runs for at}\\
\mbox{least $m$ steps}
\end{array}
\right]\\
&\ge&\Pr[ \mathcal{S} \ge1] - \Pr\left[
\begin{array}{c}
\mbox{process runs for}\\
\mbox{less than $m$ steps}
\end{array}
\right]\\
&\ge& 1- (1- p_{\mathrm{imp}})^m - \exp(-\Omega(n)).
\end{eqnarray*}
Hence, the probability that a PNE exists can be upper bounded by
\begin{eqnarray*}
&&\exp(-\Omega(n)) + \biggl(1- {1 \over16} n p(1-p)^{2n}\biggr)^m\\
&&\qquad\leq\exp(-\Omega(n)) + \exp
\bigl(-\Omega\bigl( m n p (1-p)^{2n}\bigr)\bigr) \\
&&\qquad\leq\exp\bigl(-\Omega\bigl( m n p (1-p)^{2n}\bigr)\bigr).
\end{eqnarray*}

For $p \leq1/n$ the last expression is
\[
\exp(-\Omega(n^2 p)),
\]
while for $p = g(n)/n$ where $g(n) \geq1$ the expression is
\[
\exp\bigl(-\Omega\bigl(n (1-p)^{2n}\bigr)\bigr) = \exp\bigl(-\Omega\bigl(n e^{-2 g(n)}\bigr)\bigr)=
\exp\bigl(-\Omega\bigl(e^{\log_e (n) -2 g(n)}\bigr)\bigr).
\]
This completes the proof of Theorem~\ref{thm:mediumconn}.
\end{pf*}

\subsection{Low connectivity}

\mbox{}

\begin{pf*}{Proof of Theorem~\ref{thm:lowconn}}
Note that if the graphical game is comprised of isolated edges that are
not matching pennies games, then a PNE exists. (This can be checked
easily by enumerating all best response tables for a $2 \times2$
game.) We wish to lower bound the probability of this event. To do
this, it is convenient to sample the graphical game in two stages as
follows: at the first stage we decide for each of the possible $n
\choose2$ edges whether the edge is \textit{present} (with probability
$p$) and whether it is \textit{predisposed} to be a matching pennies game
(independently with probability $1/8$); by ``predisposed'' we mean that
the edge will be set to be a matching pennies game if the edge turns
out to be isolated. At the second stage, we do the following: for an
edge that is both isolated and predisposed, we assign random payoff
tables to its endpoints conditioning on the resulting game being a
matching pennies game; for an isolated edge that is not predisposed, we
assign random payoff tables to its endpoints conditioning on the
resulting game \textit{not} being a matching pennies game; finally, for
any node that is part of a connected component with $0$ or at least $2$
edges we assign random payoff tables to the node. The probability that
there is no edge in the first stage that is both present and
predisposed is\looseness=-1
\[
(1-p/8)^{n \choose2}.
\]\looseness=0
Conditioning on this event, all present edges are not predisposed. Note
also that, when $c$ is fixed, the probability that there exists a pair
of adjacent edges is $o(1)$. It follows that the probability that all
present edges are not predisposed and no pair of edges intersect can be
lower bounded as
\[
(1-p/8)^{n \choose2} - o(1) = \biggl(1-\frac{c}{8 n^2}\biggr)^{
{n (n-1)}/{2}} - o(1).
\]
But, as explained above if all edges are isolated and none of them is a
matching pennies game a\vadjust{\goodbreak} PNE exists. Hence, the probability that a PNE
exists is at least
\[
\biggl(1-\frac{c}{8 n^2}\biggr)^{{n (n-1)}/{2}} - o(1)
\longrightarrow e^{-{c}/{16}}.
\]
\upqed\end{pf*}


%



%


\section{Deterministic graphs}

\subsection{A sufficient condition for existence of equilibria: Strong
connectivity}%
\mbox{}\vspace*{-6pt}
%
%
%
%
%
%
\begin{pf*}{Proof of Theorem~\ref{thm:expander}}
We use the same notation as in the proof of Theorem~\ref{th:
high connectivity}, except that we make the slight modification of
setting $N:=2^n-1$. Recall that $X_i$, $i=0,1,\ldots,N-1$, is the
indicator random variable of the event that the
strategy profile encoded by the number $i$ is a PNE. It is rather
straightforward (see the proof of Theorem~\ref{th:
high connectivity}) to show that
\[
\stexp[ Z ] = \stexp\Biggl[ \sum_{i=0}^{N-1}{X_i}\Biggr]=1.
\]

As in the proof of Theorem~\ref{th: high connectivity}, to
establish our result, it suffices to bound the following quantities:
\begin{eqnarray*}
b_1(G) &=& \sum_{i=0}^{N-1} \sum_{j \in B_i} \Prob[X_i=1]\Prob
[X_j=1], \\
b_2(G)&=& \sum_{i=0}^{N-1} \sum_{j \in B_i \setminus\{i\}}
\Prob[X_i=1,X_j=1],
\end{eqnarray*}
where the neighborhoods of dependence $B_i$ are defined as in Lemma
\ref{lem:neighborhood of dependence}. For $S \subseteq\{1, \ldots,
n \}$, denote by $i(S)$ the strategy profile in which the players of
the set $S$ play $1$ and the players not in $S$ play $0$. Then
writing $1(j \in B)$ for the indicator of the event that $j \in B$ we have
\begin{eqnarray*}
b_2(G) &=& \sum_{i=0}^{N-1} \sum_{j \in B_i \setminus\{i\}} \Prob
[X_i=1,X_j=1]\\
&=& \sum_{i=0}^{N-1} \sum_{j \neq i} \Prob[X_i=1,X_j=1] 1(j \in B_i)
\\
&=& N \sum_{j \neq0} \Prob[X_0=1,X_j=1] 1(j \in
B_0) \qquad\mbox{(by symmetry)}\\
&=& N \sum_{k = 1}^n \sum_{S,|S|=k} \Prob\bigl[X_0=1,X_{i(S)}=1\bigr] 1\bigl(i(S)
\in B_0\bigr).
\end{eqnarray*}
We will bound the sum above by bounding
%
%
\begin{equation} \label{eq:sum_1_ex}
N \sum_{k = 1}^{\lfloor\delta n \rfloor}
\sum_{S,|S|=k} \Prob\bigl[X_0=1,X_{i(S)}=1\bigr] 1\bigl(i(S) \in B_0\bigr)
\end{equation}
and
%
%
\begin{equation} \label{eq:sum_2_ex}
N \sum_{k = \lfloor\delta n \rfloor+1}^n \sum_{S,|S|=k} \Prob
\bigl[X_0=1,X_{i(S)}=1\bigr]
1\bigl(i(S) \in B_0\bigr)
\end{equation}
separately.

Note that if some set $S$ satisfies $|S|\le\lfloor\delta n
\rfloor$, then $|\mathcal{N}(S)| \ge\alpha|S|$ since the graph has
$(\alpha, \delta)$-expansion. Moreover, each vertex (player)
of the set $\mathcal{N}(S)$ is playing its best response to the
strategies of its neighbors in both profiles~$0$ and
$i(S)$ with probability $\frac{1}{4}$, since its
environment is different in the two profiles. On the other hand,
each player not in that set is in best response in both profiles $0$
and $i(S)$ with probability at most $\frac{1}{2}$. Hence, we can
bound~(\ref{eq:sum_1_ex}) by
\begin{eqnarray*}
&&
N \sum_{k = 1}^{\lfloor\delta n \rfloor} \sum_{S,|S|=k} \Prob
\bigl[X_0=1,X_{i(S)}=1\bigr] \\
&&\qquad\le N \sum_{k = 1}^{\lfloor\delta n \rfloor}
\sum_{S,|S|=k} \biggl( \frac{1}{2}\biggr)^{n-\alpha k} \biggl( \frac
{1}{4}\biggr)^{\alpha k}
= \sum_{k = 1}^{\lfloor\delta n \rfloor} \pmatrix{n \cr k } \biggl(
\frac{1}{2}\biggr)^{\alpha k} \\
&&\qquad< \biggl(1 + \biggl( \frac{1}{2}
\biggr)^{\alpha}\biggr)^n - 1 \leq e n^{-\varepsilon}.
\end{eqnarray*}
%
To bound the second term, notice that, if some set $S$
satisfies $|S| \ge\lfloor\delta n \rfloor+1$, then
since the graph has $(\alpha,
\delta)$-expansion $\mathcal{N}(S)\equiv V$, and, therefore,
the environment of every player is different in the two profiles $0$
and $i(S)$. Hence, $1(i(S) \in B_0) = 0$. By
combining the above we get that
\[
b_2(G) \le e n^{-\varepsilon}.
\]
It remains to bound the expression $b_1(G)$. We have
\begin{eqnarray*}
b_1(G) - 2^{-n} &=& \sum_{i=0}^{N-1} \sum_{j \in B_i \setminus\{
i\}}
\Prob[X_i=1]\Prob[X_j=1]\\
&=& \sum_{i=0}^{N-1} \sum_{j \neq i} \Prob[X_i=1] \Prob
[X_j=1] 1(j \in B_i) \\
&=& 2^{-n} \sum_{j \neq0} 1(j \in B_0) \\
&=& 2^{-n} \sum_{k = 1}^{\lfloor\delta n \rfloor} \sum
_{S,|S|=k} 1\bigl(i(s) \in B_0\bigr)\\
&&{} +
2^{-n} \sum_{k = \lfloor\delta n \rfloor+1}^n \sum_{S,|S|=k} 1\bigl(i(s)
\in B_0\bigr).
\end{eqnarray*}
The second term is zero as before. For all large enough $n$ the
first summation contains at most $2^{n/2}$ terms and is
therefore bounded by $2^{-n/2}$. It follows that
\[
b_1(G) + b_2(G) \le e n^{-\varepsilon} + 2^{-n/2}.
\]
An application of the result by Arratia, Goldstein and Gordon \cite
{Arratia89} concludes the proof of Theorem~\ref{thm:expander}.
\end{pf*}

\subsection{A sufficient condition for the nonexistence of
equilibria: Indifferent matching pennies}

In this section we provide a proof of Theorem~\ref{th:easy matching
pennies}. Recall that an edge of a graph is called \textit{d-bounded}
if both adjacent vertices have degrees smaller or equal to $d$.
Theorem~\ref{th:easy matching pennies} specifies that any graph
with many such edges is unlikely to have PNE. We proceed to the
proof of the claim.
\begin{pf*}{Proof of Theorem~\ref{th:easy matching pennies}}
Consider a $d$-bounded edge in a game connecting two players $a$
and $b$; suppose that each of these players interacts with $d-1$ (or
fewer) other players
denoted by $a_1, a_2, \ldots, a_{d-1}$ and $b_1, b_2, \ldots,
b_{d-1}$.\footnote{We allow these lists to share players.}
Recall that if $a$ and $b$ play an indifferent matching
pennies game against each other then the game has no PNE. The key
observation is that a
$d$-bounded edge is an indifferent matching pennies game with probability
at least $(\frac{1}{8})^{2^{2d-2}}=:p_{\mathrm{imp}}$---since a random
two-player game is a matching pennies game with probability
$\frac{1}{8}$, and there are at most $2^{2d-2}$ possible pure strategy
profiles for the players $a_1, a_2, \ldots, a_{d-1}$, $b_1, b_2,
\ldots, b_{d-1}$; for each of these pure strategy profiles the game
between $a$
and $b$ must be a matching pennies game.

For a collection of $m$ vertex disjoint edges, observe that the events that
each of them is an indifferent matching pennies game are independent. Hence,
the probability that the game has a PNE is upper bounded by the
probability that none of these edges is an indifferent matching pennies
game, which is upper bounded by
\[
(1-p_{\mathrm{imp}})^m \leq\exp(-m p_{\mathrm{imp}}) = \exp\bigl(-m
\bigl(\tfrac{1}{8}\bigr)^{2^{2d-2}}\bigr).\vadjust{\goodbreak}
\]

For the second claim of the theorem note that, if there are $m$ $d$-bounded
edges, then there must be at least $m/(2d)$ vertex disjoint $d$-bounded edges.

The algorithmic statement follows from the fact that we may find all
nodes with degree \mbox{$\le$}$d$ in time $O(n^2)$, and then find all edges
joining two such nodes in another $O(n^2)$ time, with the use of the
appropriate data structures; these edges are the $d$-bounded edges of
the graph. Then in time $O(m 2^{d+2})$ we can check if the endpoints of
any such edge play an indifferent matching pennies game.

The final claim of the theorem has a similar proof where now the
potential witnesses for the nonexistence of a PNE are the edges in $\E$.
\end{pf*}

Many random graphical games on deterministic graphs such as players
arranged on a line, grid, or any other bounded degree graph [with
$\omega(1)$ edges]
are special cases of the above theorem and hence
are unlikely to have PNE asymptotically.

\begin{appendix}\label{app}
\section*{Appendix: Omitted proofs}
\begin{pf*}{Proof of Lemma~\ref{lem:bounding expectations of b1 and
b2}}
We need to bound the functions $S(n,p)$ and $R(n,p)$. We begin with $S$.
\subsection*{Bounding $S$}
\label{app:high density proof}

Recall that
\[
S(n,p):=\sum_{s=1}^n \pmatrix{n \cr s} 2^{-n}
\bigl[\bigl(1+(1-p)^s\bigr)^{n-s} -
\bigl(1-(1-p)^s\bigr)^{n-s}\bigr].
\]
We split the range of the summation into four regions and
bound the sum over each region separately. We begin by choosing $\alpha
= \alpha(\varepsilon)$ as follows:
\begin{enumerate}[(ii)]
\item[(i)] if $\varepsilon\le{1790 \over105}$, we choose $\alpha=
({\varepsilon\over2+\varepsilon})^{20}$;
\item[(ii)] if $\varepsilon>{1790 \over105}$, we choose $\alpha
={\varepsilon\over2+\varepsilon}$.
\end{enumerate}
Given our choice of $\alpha=\alpha(\varepsilon)$ we define the
following regions in the range of $s$ (where---depending on $\varepsilon
$---Regions I and/or III may be empty and Region IV may have overlap
with Region II):
%
\begin{enumerate}[III.]
\item[I.] $\{ s \in\mathbb{N} | 1 \le s < \frac{\varepsilon
}{(2+\varepsilon) p}\}$;

\item[II.] $\{ s \in\mathbb{N} | \frac{\varepsilon}{(2+\varepsilon)p}
\le s < \alpha n \}$;

\item[III.] $\{ s \in\mathbb{N} | \alpha n \le s < \frac
{1}{2+\varepsilon} n\}$;

\item[IV.] $\{ s \in\mathbb{N} | \frac{1}{2+\varepsilon} n \le s <n\}$.
\end{enumerate}
We then write
\[
S(n,p) \le S_{\mathrm{I}}(n,p) + S_{\mathrm{II}}(n,p) + S_{
\mathrm{III}}(n,p) + S_{\mathrm{IV}}(n,p),
\]
where $S_{\mathrm{I}}(n,p)$ denotes the sum over region I etc., and
bound each term separately.

\subsubsection*{Region \textup{I}}
The following lemma will be useful.
\begin{lemma} \label{lem:easy_lemma}
For all $\varepsilon> 0$, $p \in(0,1)$ and $s$ such that $1 \le s <
\frac{\varepsilon}{(2+\varepsilon) p}$,
\[
(1-p)^s \le1 - \frac{(2+0.5 \varepsilon)sp}{2+ \varepsilon}.
\]
\end{lemma}
\begin{pf}
First note that, for all $k \ge1$,
%
%
\begin{equation}\label{eq:lem_easy}
\pmatrix{s \cr2k+2} p^{2k+2} \le\pmatrix{s \cr2k+1} p^{2k+1}.
\end{equation}
To verify the latter note that it is equivalent to
\[
s\le2k+1 + \frac{2k+2}{p},
\]
which is true since $s \le\frac{\varepsilon}{(2+\varepsilon) p} = \frac
{1}{({2}/{\varepsilon}+1) p}\le\frac{1}{p}$.

Using (\ref{eq:lem_easy}), it follows that
%
%
\begin{equation}\label{eq:lem_easy2}
(1-p)^s \le1 - \pmatrix{s \cr1}p+\pmatrix{s \cr2}p^2.
\end{equation}
Note finally that
\[
\frac{0.5 \varepsilon}{2+\varepsilon}sp >
\frac{s(s-1)}{2}p^2,
\]
which applied to (\ref{eq:lem_easy2}) gives
\[
(1-p)^s \le1 - \frac{(2+0.5\varepsilon)sp}{2+ \varepsilon}.
\]
\upqed\end{pf}

Assuming that Region I is nonempty and applying Lemma \ref
{lem:easy_lemma} we get
\begin{eqnarray*}
S_{\mathrm{I}}(n,p) &\le&\sum_{s < {\varepsilon}/((2+\varepsilon) p)}
\pmatrix{n \cr s} 2^{-n} \bigl(1+(1-p)^s\bigr)^{n-s} \\
&\le&\sum_{s < {\varepsilon}/((2+\varepsilon) p)} \pmatrix{n \cr s}
2^{-n} \biggl(1+1 - \frac{(2+0.5\varepsilon)sp}{2+ \varepsilon}
\biggr)^{n-s} \\
&\le&\sum_{s < {\varepsilon}/((2+\varepsilon) p)} \pmatrix{n \cr s}
2^{-s} \biggl(1 - \frac{(1+0.25\varepsilon)sp}{2+\varepsilon}
\biggr)^{n-s} \\
&\le&\sum_{s < {\varepsilon}/((2+\varepsilon) p)} \pmatrix{n \cr s} 2^{-s}
\exp\biggl(-\frac{(1+0.25\varepsilon)sp}{2+\varepsilon} (n-s) \biggr) \\
&\le&\sum_{s < {\varepsilon}/((2+\varepsilon) p)} \pmatrix{n \cr s}
2^{-s} \exp\biggl(-\frac{(1+0.25\varepsilon)sp}{2+\varepsilon}
n\biggr)\\[-2pt]
&&\hspace*{50.4pt}{}\times
\exp\biggl(\frac{(1+0.25\varepsilon)sp}{2+\varepsilon} s\biggr)\\
&\le&\sum_{s < {\varepsilon}/((2+\varepsilon) p)} \pmatrix{n \cr s} 2^{-s}
\exp\bigl(-{(1+0.25\varepsilon) \log_e (n) } s \bigr)\\[-2pt]
&&\hspace*{50.4pt}{}\times \exp
\biggl(\frac{(1+0.25\varepsilon)\varepsilon}{(2+\varepsilon)^2} s\biggr)\\[-2pt]
&\le&\sum_{s < {\varepsilon}/((2+\varepsilon) p)}{ n^s 2^{-s}
n^{-(1+0.25\varepsilon)s} \exp\biggl(\frac{1}{2}s\biggr) }\\[-2pt]
&\le&\sum_{s < {\varepsilon}/((2+\varepsilon) p)}{ \biggl(\frac
{\sqrt{e}}{2}\biggr)^{s} n^{-0.25\varepsilon s} }\\[-2pt]
&\le&\sum_{s < {\varepsilon}/((2+\varepsilon) p)}{ \biggl(\frac
{\sqrt{e}}{2}\biggr)^{s} n^{-0.25\varepsilon} }\\[-2pt]
&\le& n^{-0.25\varepsilon} \sum_{s < {2 \varepsilon}/({(2+\varepsilon)
p})}{ \biggl(\frac{\sqrt{e}}{2}\biggr)^{s}}\\[-2pt]
&=& O(n^{-0.25\varepsilon})\qquad \biggl(\mbox{since $\dfrac
{\sqrt{e}}{2} < 1$}\biggr).
\end{eqnarray*}

\subsubsection*{Region \textup{II}}
We have
%
%
\begin{eqnarray} \label{eq:prefinal region II}
S_{\mathrm{II}}(n,p) &\le&\sum_{{ \varepsilon}/({(2+\varepsilon)p})
\le
s < \alpha n} \pmatrix{n \cr s} 2^{-n} \bigl(1+(1-p)^s\bigr)^n \nonumber\\[-2pt]
&\le&\sum_{{\varepsilon}/({(2+\varepsilon)p}) \le s < \alpha n}
\pmatrix{n
\cr s} 2^{-n} (1+e^{-ps})^n \nonumber\\[-2pt]
&\le&\sum_{{ \varepsilon}/({(2+\varepsilon)p}) \le s < \alpha n}
\pmatrix{n
\cr\alpha n} 2^{-n} \bigl(1+e^{-p{\varepsilon}/({(2+\varepsilon
)p})}\bigr)^n \nonumber\\[-9pt]\\[-9pt]
&\le&\alpha n \pmatrix{n \cr\alpha n} \biggl( \frac{1+e^{-{
\varepsilon}/({2+\varepsilon})}}{2}\biggr)^n \nonumber\\[-2pt]
&\le&\alpha n 2^{nH(\alpha)} (n+1)\biggl( \frac{1+e^{-{\varepsilon
}/({2+\varepsilon})}}{2}\biggr)^n \nonumber\\[-2pt]
&\le&\alpha n (n+1) \biggl( 2^{H(\alpha)} \cdot\frac{1+e^{-{
\varepsilon}/({2+\varepsilon})}}{2}\biggr)^n.\nonumber
\end{eqnarray}
In the above derivation $H(\cdot)$ represents the entropy function,
and for the second to last derivation we used the fact that
%
%
\begin{equation}\label{eq:entropy bound}
\pmatrix{n \cr k} \le(n+1)2^{nH({k}/{n})}.
\end{equation}
%
Our definition of the function $\alpha=\alpha(\varepsilon)$ guarantees
that when $\varepsilon\le{1790 \over105}$
\[
\biggl( 2^{H(\alpha)} \cdot\frac{1+e^{-{ \varepsilon}/({2+\varepsilon
})}}{2}\biggr) \le0.999,
\]
while when $\varepsilon>{1790 \over105}$
\[
\biggl( 2^{H(\alpha)} \cdot\frac{1+e^{-{ \varepsilon}/({2+\varepsilon
})}}{2}\biggr) \le0.99.
\]
Using the above and (\ref{eq:prefinal region II}) we obtain
%
%
\begin{equation}
S_{\mathrm{II}}(n,p) = \exp(-\Omega(n)).
\end{equation}

\subsubsection*{Region \textup{III}}
Let us assume that the region is nonempty. We show that each positive
term in the summation $S_{\mathrm{III}}(n,p)$ is exponentially small.
Since there are $O(n)$ terms in the summation, it follows then that
$S_{\mathrm{III}}(n,p)$ is exponentially small:
%
%
\begin{eqnarray} \label{eq:lalalalal}
&&\pmatrix{n \cr s} 2^{-n} \bigl(1+(1-p)^s\bigr)^n \nonumber\\
&&\qquad\le\pmatrix{n \cr s} 2^{-n}
(1+e^{-ps})^{n}
\le\pmatrix{n \cr s} 2^{-n} (1+e^{-p \alpha n})^{n} \\
&&\qquad\le\pmatrix{n \cr s} 2^{-n} \bigl(1+e^{-(2+\varepsilon) \alpha\log_e
(n)}\bigr)^{n} \nonumber\\
&&\qquad= \pmatrix{n \cr s} 2^{-n} \biggl(1+\frac{1}{n^{(2+\varepsilon) \alpha
}}\biggr)^{n} \nonumber\\
&&\qquad= \pmatrix{n \cr s} 2^{-n} \biggl(1+\frac{1}{n^{(2+\varepsilon) \alpha
}}\biggr)^{n^{(2+\varepsilon) \alpha} n^{1-(2+\varepsilon) \alpha}}
\nonumber\\
&&\qquad\le\pmatrix{n \cr s} 2^{-n} e^{n^{1-(2+\varepsilon) \alpha}} \nonumber
\\
&&\qquad\le(n+1)2^{n H ({s}/{n})} 2^{-n}
e^{n^{1-(2+\varepsilon) \alpha}} \nonumber\\
&&\qquad\le(n+1)2^{n H ({1}/({2+\varepsilon}))} 2^{-n}
e^{n^{1-(2+\varepsilon) \alpha}} \nonumber\\
&&\qquad = (n+1)2^{n (H ({1}/({2+\varepsilon}))-1)}
e^{n^{1-(2+\varepsilon) \alpha}},\nonumber
\end{eqnarray}
where in the third-to-last line of the derivation we employed the bound
of~(\ref{eq:entropy bound}). Notice that the right-hand\vadjust{\goodbreak} side of (\ref
{eq:lalalalal}), seen as a function of $\varepsilon>0$ and $\alpha>0$, is
decreasing in both. Since $\varepsilon>c$, our choice of $\alpha=\alpha
(\varepsilon)$ implies that $\alpha>( {c \over c+2} )^{20}$.
Hence, we can bound the right-hand side of (\ref{eq:lalalalal}) as follows:
\begin{eqnarray*}
&&(n+1)2^{- n (1-H ({1}/({2+c})))}
e^{n^{1-(2+c)( {c / (c+2)} )^{20}}} \\
&&\qquad= \exp(-\Omega(n)),
\end{eqnarray*}
where we used the fact that $c$ is a constant, and therefore the factor
$e^{n^{1-(2+c)( c / (c+2) )^{20}}}$ is sub-exponential
in $n$, while the factor $2^{- n (1-H ({1}/({2+c})
))}$ is exponentially small in $n$.

\subsubsection*{Region \textup{IV}}
Note that, if $x k \le1$, then by the mean value theorem
\begin{eqnarray*}
(1+ x)^k - (1- x)^k &\leq&2x \max_{1-1/k \leq y \leq1+1/k} k y^{k-1} \\
&=&
2k x (1+1/k)^{k-1} \\
&\leq&2 e k x.
\end{eqnarray*}
%

We can apply this for $k = n-s$ and $x = (1-p)^s$ since
\begin{eqnarray*}
(n-s)(1-p)^s &\le&(n-s)e^{-ps} \\
&\le&(n-s)e^{- {(2+\varepsilon)\log_e
(n)}/{n}({n}/({2 + \varepsilon}))} \\
&\le&\frac{n-s}{n} \\
&\le&1.
\end{eqnarray*}

Hence, $S_{\mathrm{IV}}(n,p)$ is bounded as follows:
\begin{eqnarray*}
S_{\mathrm{IV}}(n,p) &\le& \sum_{{n}/({2+\varepsilon}) \le s \le
n}\pmatrix{n
\cr s} 2^{-n} 2 e (n-s)(1-p)^s\\
&\le& 2 e \cdot2^{-n} \cdot n \sum_{{n}/({2+\varepsilon})\le s \le
n}\pmatrix{n \cr s} (1-p)^s\\
&\le& 2 e \cdot2^{-n} \cdot n \bigl(1+(1-p)\bigr)^n\\
&\le& 2 e n \biggl(1-\frac{p}{2}\biggr)^n\\
&\le& 2 e n e^{-{p}/{2}n}\\
&\le& 2 e n e^{-{(2+\varepsilon)\log_e (n)}/({2 n})n}\\
&\le& 2 e n n^{-({2+\varepsilon})/{2}}\\
&\le& 2 e n^{-{\varepsilon}/{2}}.
\end{eqnarray*}

\subsubsection*{Putting everything together}
Combining the above we get that
\[
S(n,p) \leq O(n^{-\eps/4}) + \exp(-\Omega(n)).
\]

\subsection*{Bounding $R$}
Observe that
\[
R(n,p) =2^{-n}+\sum_{s=1}^n \pmatrix{n \cr s} 2^{-n} \min\bigl(1, n (1-p)^{s-1}\bigr).
\]
We bound $R$ as follows:
\begin{eqnarray*}
&&
R(n,p) -2^{-n} \\[-2pt]
&&\qquad\leq\sum_{s=1}^n \pmatrix{n \cr s} 2^{-n} \min\bigl(1, n
\exp\bigl(-p(s-1)\bigr)\bigr) \\[-2pt]
&&\qquad \leq2^{-n} \sum_{1 \leq s \leq{n(3+\eps)/{(6+3\eps)}}} \pmatrix{ n \cr s}\\
&&\qquad\quad{} +
2^{-n} \sum_{s>{n(3+\eps)/{(6+3\eps)}}} \pmatrix{n \cr s} n
\exp\bigl(-p(s-1)\bigr) \\[-2pt]
&&\qquad
\leq2^{-n} \sum_{1 \leq s \leq{n(3+\eps)/{(6+3\eps)}}} {
(n+1)2^{nH(s/n)}}\\[-2pt]
&&\qquad\quad{}+2^{-n} \sum_{s > {n(3+\eps)/{(6+3\eps)}}} \pmatrix{n \cr s} n \exp\bigl(-p(s-1)\bigr) \\[-2pt]
&&\qquad \leq n (n+1) 2^{-n} {2^{nH({(3+\eps)/({6+3\eps})}
)}} \\[-2pt]
&&\qquad\quad{}+2^{-n} \sum_{s > {n(3+\eps)/{(6+3\eps)}}} \pmatrix{n
\cr
s} n \exp\bigl(-p(s-1)\bigr)\\[-2pt]
&&\qquad \leq\exp(-\Omega(n)) +2^{-n} \sum_{s >{n(3+\eps)/{(6+3\eps)}}} \pmatrix{n \cr s} n \exp\bigl(-p(s-1)\bigr),
\end{eqnarray*}
where in the last line of the derivation we used that $\varepsilon>c>0$
for some absolute constant $c$.
To bound the last sum we observe that when $s > {n(3+\eps) \over{6+3\eps
}}$ we
have
\begin{eqnarray*}
n \exp\bigl(-p(s-1)\bigr) &\leq& n \exp\biggl(-\frac{(2+\eps) \log_e
(n)}{n}\biggl( {n(3+\eps) \over{6+3\eps}}-1\biggr)\biggr)\\
&\leq& n \cdot n^{-{(2+\eps)(3+\eps)/{(6+3\eps)}}} \cdot\exp
\biggl(\frac{(2+\eps) \log_e (n)}{n}\biggr)\\
&\leq& n^{-\eps/3} \cdot n^{2/n} \cdot n^{\eps/n} =O(n^{-\eps/4}).
\end{eqnarray*}
Using this bound and the fact $\sum_{s=0}^n {n \choose s}=2^n$
concludes the proof.\vadjust{\goodbreak}
\end{pf*}
%
\end{appendix}

\section*{Acknowledgments}
We thank Martin Dyer for pointing out an error in a previous
formulation of Theorem~\ref{th:easy matching pennies}. We also thank
the anonymous referee for comments that helped improve the presentation
of this work.

%

%
\printaddresses

\end{document}